\documentclass[11pt]{amsart}
\usepackage{graphicx}
\usepackage{latexsym}
\usepackage{amssymb}
\usepackage{amsfonts}
\usepackage{amsmath}
\usepackage{amscd}

\title{A Statement of the Fundamental Lemma}
\date{30 December 2003}

\author{Thomas C. Hales}

\address{Department of Mathematics, University of Pittsburgh,
Pittsburgh, PA 15260} \email{hales@pitt.edu}

\issueinfo{00}% volume number
{0}         % issue number
{0}         % month
{}      % year
\copyrightinfo{2003}% copyright year
{Hales}% copyright holder

\newtheorem{theorem}{Theorem}[section]

\newtheorem{conjecture}[theorem]{Conjecture}

\newtheorem{assumption}[theorem]{Assumption}

\theoremstyle{definition}
\newtheorem{definition}[theorem]{Definition}
\newtheorem{example}[theorem]{Example}
\theoremstyle{remark}
\newtheorem{remark}[theorem]{Remark}
\numberwithin{equation}{theorem}

\def\tq{{\ \vert\ }}

\def\calG{\mathcal{G}}

\newcommand{\ring}[1]{\mathbb{#1}}

\def\op#1{\operatorname{#1}}

\def\trans{{{}^t}}

\begin{document}

\begin{abstract}
These notes give a statement of the fundamental lemma, which is a
conjectural identity between $p$-adic integrals.
\end{abstract}

\thanks{I thank R. Kottwitz and S. DeBacker for many
helpful comments.}

\maketitle

\def\today{\ifcase\month\or
    January\or February\or March\or April\or May\or June\or
    July\or August\or September\or October\or November\or December\fi
    \space\number\day, \number\year}

%\begin{document}
\maketitle

\section{Introduction}\label{section: introduction}

\noindent{\bf Notation.\/}  Let $F$ be a $p$-adic field, given
either as a finite field extension of $\ring{Q}_p$, or as the
field $F=\ring{F}_q((t))$. Let $\ring{F}_q$ (a finite field with
$q$ elements and characteristic $p$) be the residue field of $F$.
Let $\bar F$ be a fixed algebraic closure of $F$.  Let $F^{un}$ be
the maximal unramified extension of $F$ in $\bar F$.  For
simplicity, we also assume that the characteristic of $F$ is not
$2$.

The fundamental lemma pertains to groups that satisfy a series of
hypotheses. Here is the first.

\begin{assumption}\label{ass:reductive}
$G$ is a connected reductive linear algebraic group that is
defined over $F$.
\end{assumption}

The following examples give the $F$-points of three different
families of connected reductive linear algebraic groups:
orthogonal, symplectic, and unitary groups.

\begin{example}  Let $M(n,F)$ be the algebra of $n$ by $n$ matrices with
coefficients in $F$.  Let $J\in M(n,F)$ be a symmetric matrix with
nonzero determinant. The special orthogonal group with respect to
the matrix $J$ is
    $$
    \op{SO}(n,J,F) = \{ X \in M(n,F) \tq \trans X J X = J,\quad
    \op{det}(X)=1\}.
    $$
\end{example}

\begin{example} Let $J\in M(n,F)$, with $n=2k$, be a skew-symmetric matrix
$\trans J = -J$ with nonzero determinant.  The symplectic group
with respect to $J$ is defined in a similar manner:
    $$
    \op{Sp}(2k,J,F) = \{ X \in M(2k,F) \tq \trans X J X = J
    \}.
    $$
\end{example}

\begin{example} Let $E/F$ be a separable quadratic extension.  Let $\bar x$
be the Galois conjugate of $x\in E$ with respect to the nontrivial
automorphism of $E$ fixing $F$. For any $A\in M(n,E)$, let $\bar
A$ be the matrix obtained by taking the Galois conjugate of each
coefficient of $A$.  Let $J\in M(n,E)$ satisfy $\trans {\bar J} =
J$ and have a nonzero determinant.  The unitary group with respect
to $J$ and $E/F$ is
    $$
    U(n,J,F) = \{X \in M(n,E) \tq \trans {\bar X} J X = J\}.
    $$
\end{example}

The algebraic groups $SO(n,J)$, $Sp(2k,J)$, and $U(n,J)$ satisfy
Assumption~\ref{ass:reductive}.

\begin{assumption}\label{ass:unram} $G$ splits over an unramified field
extension.
\end{assumption}

That is, there is an unramified extension $F_1/F$ such that
$G\times_F F_1$ is split.

\begin{example}\label{ex:J}  In the first two examples above (orthogonal and
symplectic), if we take $J$ to have the special form
    \begin{equation}
    J = \begin{pmatrix}0&0&*\\0&*&0\\{*}&0&0\end{pmatrix}
    \label{eqn:cross}
    \end{equation}
(that is, nonzero entries from $F$ along the cross-diagonal and
zeros elsewhere), then $G$ splits over $F$. In the third example
(unitary), if $J$ has this same form  and if $E/F$ is unramified,
then the unitary group splits over the unramified extension $E$ of
$F$.
\end{example}

\begin{assumption}\label{ass:quasi-split} $G$ is quasi-split.
\end{assumption}

This means that there is an $F$-subgroup $B\subset G$ such that
$B\times_F\bar F$ is a Borel subgroup of $G\times_F\bar F$.

\begin{example} In all three cases (orthogonal, symplectic, and
unitary), if $J$ has the cross-diagonal form \ref{eqn:cross}, then
$G$ is quasi-split.  In fact, we can take the points of $B$ to be
the set of upper triangular matrices in $G(F)$.
\end{example}

\begin{assumption}  $K$ is a hyperspecial maximal
compact subgroup of $G(F)$, in the sense of Definition
\ref{def:hyper}.
\end{assumption}

\begin{example} Let $O_F$ be the ring of integers of $F$ and let $K = GL(n,O_F)$.
This is a hyperspecial maximal compact subgroup of $GL(n,F)$.
\end{example}

\begin{definition}\label{def:hyper}
$K$ is {\it hyperspecial} if there exists $\calG$ such that the
following conditions are satisfied.
    \begin{itemize}
    \item $\calG$ is a smooth group scheme over $O_F$,
    \item $G = \calG\times _{O_F} F$,
    \item $\calG\times_{O_F} \ring{F}_q$ is reductive,
        \item $K = \calG(O_F)$.
    \end{itemize}
\end{definition}

\begin{example} In all three examples (orthogonal, symplectic, and
unitary), take $G$ to have the form of Example~\ref{ex:J}. Assume
that each cross-diagonal entry is a {\it unit} in the ring of
integers. Assume further that the residual characteristic is not
$2$.  Then the equations
    $$
    \trans X J X = J\quad (\text{or in the unitary case } \trans{\bar X} J X = J)
    $$
define a group scheme $\calG$ over $O_F$, and $\calG(O_F)$ is
hyperspecial.
\end{example}

\section{Classification of Unramified Reductive Groups}

\begin{definition} If $G$ is quasi-split and splits over an unramified extension
(that is, if $G$ satisfies Assumptions~\ref{ass:unram} and
\ref{ass:quasi-split}), then $G$ is said to be an {\it unramified
reductive group}.
\end{definition}

Let $G$ be an unramified reductive group.  It is classified by
data (called root data)
    $$(X^*,X_*,\Phi,\Phi^{\vee},\sigma).$$
The data is as follows:
    \begin{itemize}
    \item $X^*$ is the character group of a Cartan
subgroup of $G$.
    \item
$X_*$ is the cocharacter group of the Cartan subgroup.
    \item
$\Phi\subset X^*$ is the set of roots.
    \item
$\Phi^{\vee} \subset X_*$ is the set of coroots.
    \item
$\sigma$ is an automorphism of finite order of $X^*$ sending a set
of simple roots in $\Phi$ to itself.

$\sigma$ is obtained from the action on the character group
induced from the Frobenius automorphism of $\op{Gal}(F^{un}/F)$ on
the maximally split Cartan subgroup in $G$.
    \end{itemize}

The first four elements $(X^*,X_*,\Phi,\Phi^{\vee})$ classify
split reductive groups $G$ over $F$.  For such groups $\sigma=1$.

\section{Endoscopic Groups}

$H$ is an unramified endoscopic group of $G$ if it is an
unramified reductive group over $F$ whose classifying data has the
form
    $$
    (X^*,X_*,\Phi_H,\Phi_H^{\vee},\sigma_H).
    $$
The first two entries are the same for $G$ as for $H$. To
distinguish the data for $H$ from that for $G$, we add subscripts
$H$ or $G$, as needed. The data for $H$ is subject to the
constraints that there exists an element $s \in
\op{Hom}(X_*,\ring{C}^\times)$ and a Weyl group element $w\in
W(\Phi_G)$ such that
    \begin{itemize}
        \item$\Phi_H^{\vee}= \{\alpha\in\Phi_G^{\vee} \tq
        s(\alpha)=1\}$,
        \item $\sigma_H = w\circ \sigma_G$, and
        \item $\sigma_H(s) = s$.
    \end{itemize}

\subsection{Endoscopic groups for $SL(2)$}

As an example, we determine the unramified endoscopic groups of $G
= SL(2)$. The character group $X^*$ can be identified with
$\ring{Z}$, where $n\in\ring{Z}$ is identified with the character
on the diagonal torus given by
    $$
    \begin{pmatrix} t&0\\0&t^{-1} \end{pmatrix} \mapsto t^n.
    $$
The set $\Phi$ can be identified with the subset $\{\pm 2\}$ of
$\ring{Z}$:
    $$
    \begin{pmatrix} t&0\\0&t^{-1} \end{pmatrix} \mapsto t^{\pm 2}.
    $$
The cocharacter group $X_*$ is also identified with $\ring{Z}$,
where $n\in \ring{Z}$ is identified with
    $$
    t\mapsto \begin{pmatrix} t^n&0\\0&t^{-n}\end{pmatrix}.
    $$
Under this identification $\Phi^{\vee} = \{\pm 1\}$.  Since the
group is split, $\sigma = 1$.

We get an unramified endoscopic group by selecting
$s\in\op{Hom}(X_*,\ring{C}^\times)\cong \ring{C}^\times$ and $w\in
W(\Phi)$.
    \begin{equation}
    \begin{array}{lll}
    \Phi_H^{\vee} = \{\alpha\tq s(\alpha)=1\} &= \{n\in\{\pm1\}\tq
    s^n = 1\} \\&=
    \text{if}~~(s = 1)~~\text{then}~~\Phi_G^{\vee}~~
    \text{else}~~\emptyset.
        \end{array}
    \end{equation}

We consider two cases, according as $w$ is nontrivial or trivial.
If $w$ is the nontrivial reflection, then $\sigma_H = w$ acts by
negation on $\ring{Z}$. Thus,
    $$
    (\sigma_H(s) = s)~~ \implies~~( s^{-1} = s)~~ \implies (s = \pm 1).
    $$
If $s=1$, then $\sigma_H$ does not fix a set of simple roots as
required.  So $s=-1$ and $\Phi_H^{\vee} = \emptyset$. Thus, the
root data of $H$ is
    $$
    (\ring{Z},\ring{Z},\emptyset,\emptyset,w)
    $$
This determines $H$ up to isomorphism as $H=U_E(1)$, a
$1$-dimensional torus split by an unramified quadratic extension
$E/F$.

If $w$ is trivial, then there are two further cases, according as
$\Phi_H$ is empty or not:
    \begin{itemize}
        \item The endoscopic group $\ring{G}_m$ has root data
        $$
        (\ring{Z},\ring{Z},\emptyset,\emptyset,1).
        $$
    \item The endoscopic group $H=SL(2)$ has root data
        $$
        (\ring{Z},\ring{Z},\{\pm2\},\{\pm1\},1).
        $$
    \end{itemize}

In summary, the three unramified endoscopic groups of $SL(2)$ are
$U_E(1)$, $\ring{G}_m$, and $SL(2)$ itself.

\subsection{Endoscopic groups for $PGL(2)$}

As a second complete example, we determine the endoscopic groups
of $PGL(2)$.   The group $PGL(2)$ is dual to $SL(2)$ in the sense
that the coroots of one group can be identified with the roots of
the other group.  The root data for $PGL(2)$ is
    $$
    (\ring{Z},\ring{Z},\{\pm1\},\{\pm2\},1).
    $$

When the Weyl group element is trivial, then the calculation is
almost identical to the calculation for $SL(2)$.  We find that
there are again two cases, according as  $\Phi_H$ is empty or not:
    \begin{itemize}
   \item The endoscopic group $\ring{G}_m$ has root data
        $$
        (\ring{Z},\ring{Z},\emptyset,\emptyset,1).
        $$
    \item The endoscopic group $H=PGL(2)$ has root data
        $$
        (\ring{Z},\ring{Z},\{\pm1\},\{\pm2\},1).
        $$
    \end{itemize}

When the Weyl group element $w$ is nontrivial, then
$s\in\{\pm1\}$, as in the $SL(2)$ calculation.
    \begin{equation}
    \Phi_H^{\vee} = \{\alpha\tq s(\alpha)=1\} =
        \{n  \in\{\pm2\}\tq s^n = 1\} = \Phi_G^{\vee}.
    \end{equation}
From this, we see that picking $w$ to be nontrivial is
incompatible with the requirement that $\sigma_H=w$ must fix a set
of simple roots.  Thus, there are no endoscopic groups with $w$
nontrivial.

In summary, the two endoscopic groups of $PGL(2)$ are $\ring{G}_m$
and $PGL(2)$ itself.

\subsection{Elliptic Endoscopic groups}

\begin{definition}  An unramified endoscopic group $H$ is said to be {\it elliptic},  if
    $$
    (\ring{R}\Phi_G)^{W(\Phi_H)\rtimes\langle \sigma_H\rangle} =
    (0).
    $$
That is, the span of the set of roots of $G$ has no invariant
vectors under the Weyl group of $H$ and the automorphism
$\sigma_H$.
\end{definition}

The origin of the term {\it elliptic} is the following.  We will
see below that each Cartan subgroup of $H$ is isomorphic to a
Cartan subgroup of $G$.  (Here and elsewhere, when we speak of an
isomorphic between algebraic groups defined over $F$, we mean an
isomorphism over $F$.)  The condition on $H$ for it to be elliptic
is precisely the condition that is needed for some Cartan subgroup
of $H$ to be isomorphic to an {\it elliptic} Cartan subgroup of
$G$.

\begin{example}  We calculate the elliptic unramified endoscopic subgroups of
$SL(2)$.  We may identify $\ring{R}\Phi$ with $\ring{R}\{\pm2\}$
and hence with $\ring{R}$.  An unramified endoscopic group is
elliptic precisely when $W(\Phi_H)$ or $\langle \sigma_H\rangle$
contains the nontrivial reflection $x\mapsto -x$.  When $H=SL(2)$,
the Weyl group contains the nontrivial reflection.  When
$H=U_E(1)$, the element $\sigma_H$ is the nontrivial reflection.
But when $H=\ring{G}_m$, both $W(\Phi_H)$ and $\langle
\sigma_H\rangle$ are trivial.  Thus, $H=SL(2)$ and $H=U_E(1)$ are
elliptic, but $H=\ring{G}_m$ is not.
\end{example}

\subsection{An exercise: elliptic endoscopic groups of unitary groups}

This exercise is a calculation of the elliptic unramified
endoscopic groups of $U(n,J)$.   We assume that $J$ is a
cross-diagonal matrix with units along the cross-diagonal as in
Section~\ref{eqn:cross}. We give a few facts about the endoscopic
groups of $U(n,J)$ and leave it as an exercise to fill in the
details.

    Let $T = \{\op{diag}(t_1,\ldots,t_n)\}$ be the group of
    diagonal $n$ by $n$ matrices.
    The character group $X^*$ can be identified with $\ring{Z}^n$
in such a way that the character
    $$
    \op{diag}(t_1,\ldots,t_n) \mapsto t_1^{k_1}\cdots t_n^{k_n}
    $$
is identified with $(k_1,\ldots,k_n)$.

The cocharacter group can be identified with $\ring{Z}^n$ in such
a way that the cocharacter
    $$
    t\mapsto \op{diag}(t^{k_1},\ldots,t^{k_n})
    $$
is identified with $(k_1,\ldots,k_n)$.

Let $e_i$ be the basis vector of $\ring{Z}^n$ whose $j$-th entry
is Kronecker $\delta_{ij}$.  The set of roots can be identified
with
    $$
    \Phi = \{e_i-e_j \tq i\ne j\}.
    $$
The set of coroots $\Phi^{\vee}$ can be identified with the set of
roots $\Phi$ under the isomorphism $X_* \cong \ring{Z}^n \cong
X^*$.

We may identify  $\op{Hom}(X_*,\ring{C}^\times)$ with
$\op{Hom}(\ring{Z}^n,\ring{C}^\times)= (\ring{C}^\times)^n$. Thus,
we take the element $s$ in the definition of endoscopic group to
have the form $s=(s_1,\ldots,s_n)\in(\ring{C}^\times)^n$.  The
element $\sigma = \sigma_G$ acts on characters and cocharacters by
    $$
    \sigma(k_1,\ldots,k_n) = (-k_n,\ldots,-k_1).
    $$

Let $I = \{1,\ldots,n\}$.  Show that if $H$ is an elliptic
unramified endoscopic group, then there is a partition
    $$I = I_1\coprod I_2$$
with $s_i = 1$ for $i\in I_1$ and $s_i = -1$ otherwise.  The
elliptic endoscopic group is a product of two smaller unitary
groups $H = U(n_1)\times U(n_2)$, where $n_i = \#I_i$, for
$i=1,2$.

\section{Cartan subgroups}

All unramified reductive groups are classified by their root data.
This includes the classification of unramified tori $T$ as a
special case (in this case, the set of roots and the set of
coroots are empty):
    $$
    (X^*(T),X_*(T),\emptyset,\emptyset,\sigma).
    $$
We can extend this classification to ramified tori.   If $T$ is
any torus over $F$, it is classified by
    $$
    (X^*(T),X_*(T),\rho),
    $$
where $\rho$ is now allowed to be any homomorphism
    $$\rho:\op{Gal}(\bar F/F)\to \op{Aut}(X^*(T))$$
with finite image.

A basic fact is that $T$ embeds over $F$ as a Cartan subgroup in a
given unramified reductive group $G$ if and only if the following
two conditions hold.
    \begin{itemize}
    \item The image of $\rho$ in $\op{Aut}(X^*(T))$ is contained in
        $W(\Phi_G)\rtimes\langle\sigma_G\rangle$.
    \item
    There is a commutative diagram:
    \end{itemize}

\begin{equation*}
    \begin{CD}
        \op{Gal}(\bar F/F)@>>>\op{Gal}(F^{un}/F)\\
        @V\rho VV   @VV{\op{Frob}\mapsto \sigma_G}V\\
        W(\Phi_G)\rtimes\langle\sigma_G\rangle
        @>>w\rtimes\tau\mapsto\tau>
        \langle\sigma_G\rangle.
    \end{CD}
\end{equation*}

It follows that every Cartan subgroup $T_H$ of $H$ is isomorphic
over $F$  with a Cartan subgroup $T_G$ of $G$. (To check this,
simply observe that these two conditions are more restrictive for
$H$ than the corresponding conditions for $G$.)  The isomorphism
can be chosen to induce an isomorphism of Galois modules between
the character group (and cocharacter group) of $T_H$ and that of
$T_G$.

We say that a semisimple element in a reductive group is {\it
strongly regular}, if its centralizer is a Cartan subgroup.  If
$\gamma\in H(F)$ is strongly regular semisimple, then its
centralizer $T_H$ is isomorphic to some $T_G\subset G$. Let
$\gamma_0\in T_G(F)\subset G(F)$ be the element in $G(F)$
corresponding to $\gamma\in T_H(F)\subset H(F)$, under this
isomorphism.

\begin{remark} The element $\gamma_0$ is not uniquely determined
by $\gamma$.  The Cartan subgroup $T_G$ can always be replaced
with a conjugate $g^{-1}\,T_G\,g$, $g\in G(F)$, without altering
the root data. However, the non-uniqueness runs deeper than this.
An example will be worked in Section~\ref{sec:bp} to show how to
deal with the problem of non-uniqueness.  Non-uniqueness of
$\gamma_0$ is related to stable conjugacy, which is our next
topic.
\end{remark}

\section{Stable Conjugacy}

\begin{definition}
Let $\delta$ and $\delta'$ be strongly regular semisimple elements
in $G(F)$. They are {\it conjugate} if $g^{-1}\delta g = \delta'$
for some $g\in G(F)$.  They are {\it stably conjugate} if
$g^{-1}\delta g = \delta'$ for some $g\in G(\bar F)$.
\end{definition}

\begin{example}  Let $G= SL(2)$ and $F=\ring{Q}_p$.  Assume that
$p\ne 2$ and that $u$ is a unit that is not a square in
$\ring{Q}_p$. Let $\epsilon = \sqrt{u}$ in an unramified quadratic
extension of $\ring{Q}_p$. We have the matrix calculation
    $$
    \begin{pmatrix} 1+p&1\\2p+p^2&1+p\end{pmatrix}
    \begin{pmatrix} \epsilon&0\\0&\epsilon^{-1}\end{pmatrix}=
    \begin{pmatrix} \epsilon&0\\0&\epsilon^{-1} \end{pmatrix}
    \begin{pmatrix} 1+p&u^{-1}\\(2p+p^2)u&1+p\end{pmatrix}.
    $$
This matrix calculation shows that the matrices
    \begin{equation}
    \begin{pmatrix}1+p&1\\2p+p^2&1+p\end{pmatrix}\text{ and }
    \begin{pmatrix}1+p&u^{-1}\\(2p+p^2)u&1+p\end{pmatrix}
    \label{eqn:2matrix}
    \end{equation}
of $SL(2,\ring{Q}_p)$ are stably conjugate. The diagonal matrix
that conjugates one to the other has coefficients that lie in a
quadratic extension. A short calculation shows that the matrices
\ref{eqn:2matrix} are not conjugate by a matrix of
$SL(2,\ring{Q}_p)$.
\end{example}

\subsection{Cocycles}

Let $\gamma_0$ and $\gamma'$ be stably conjugate strongly regular
semisimple elements of $G(F)$.  We view $\gamma_0$ as a fixed base
point and $\gamma'$ as variable.  If $\tau\in \op{Gal}(\bar F/F)$,
then
    \begin{equation}
    \begin{array}{lll}
        g^{-1}\gamma_0 g &= \gamma', (\text{with } g\in G(\bar F), \gamma_0,\gamma'\in G(F))\\
        \tau(g)^{-1}\tau(\gamma_0) \tau (g) &= \tau(\gamma'),\\
        \tau(g)^{-1}\gamma_0\tau(g) &= g^{-1}\gamma_0 g,\\
        \gamma_0 \left( \tau(g)g^{-1}\right) &=
        \left(\tau(g)g^{-1}\right) \gamma_0,\\
        \gamma_0 a_\tau = a_\tau \gamma_0, \text{ with } a_\tau  =
        \tau{(g)}g^{-1}.
    \end{array}
    \end{equation}
The element $a_\tau$ centralizes $\gamma_0$ and hence gives an
element of the centralizer $T$.  Viewed as a function of $\tau\in
\op{Gal}(\bar F/F)$, $a_\tau$ satisfies the cocycle relation
    $$\tau_1(a_{\tau_2})a_{\tau_1} = a_{\tau_1\tau_2}.$$
It is continuous in the sense that there exists a field extension
$F_1/F$ for which $a_\tau=1$, for all $\tau\in\op{Gal}(\bar
F/F_1)$.   Thus, $a_\tau$ gives a class in
    $$
    H^1(\op{Gal}(\bar F/F),T(\bar F)),
    $$
which is defined to be the group of all continuous cocycles with
values in $T$, modulo the subgroup of all continuous cocycles of
the form
    $$
    b_\tau = \tau(t)t^{-1},
    $$
for some $t\in T(\bar F)$.

A general calculation of the group $H^1(\op{Gal}(\bar F/F),T)$ is
achieved by the {\it Tate-Nakayama} isomorphism.  Let $F_1/F$ be a
Galois extension that splits the Cartan subgroup $T$.

\begin{theorem} (Tate-Nakayama isomorphism \cite{SLF})
The group $H^1(\op{Gal}(\bar F/F),T)$
is isomorphic to the quotient of the group
    $$
    \{u \in X_*\tq \sum_{\tau\in\op{Gal}(F_1/F)} \tau u = 0\}
    $$
by the subgroup generated by the set
    $$
    \{u \in X_* \tq \exists \tau\in \op{Gal}(F_1/F)~~ \exists v \in
    X_*.~~ u = \tau v - v\}.
    $$
\end{theorem}

\begin{example}  Let $T = U_E(1)$ (the torus that made an
appearance earlier as an endoscopic group of $SL(2)$).  As was
shown above, the group of cocharacters can be identified with
$\ring{Z}$.  The splitting field of $T$ is the quadratic extension
field $E$.  The nontrivial element $\tau\in\op{Gal}(E/F)$ acts by
reflection on $X_*\cong\ring{Z}$:  $\tau(u) = -u$.  By the
Tate-Nakayama isomorphism, the  group $H^1(\op{Gal}(\bar
F/F),U_E(1))$ is isomorphic to
    $$
    %\begin{array}{lll}
    \{u\in\ring{Z} \tq u+\tau u = 0\}/\{u\in\ring{Z}\tq \exists v.~u = \tau v -
    v\}=\ring{Z}/2\ring{Z}.
    %\end{array}
    $$
\end{example}

Let $H$ be an unramified endoscopic group of $G$.   Suppose that
$T_H$ is a Cartan subgroup of $H$. Let $T_G$ be an isomorphic
Cartan subgroup in $G$.  The data defining $H$ includes the
existence of an element $s\in \op{Hom}(X_*,\ring{C}^\times)$; that
is, a character of the abelian group $X_*$.  Fix one such
character $s$.  We can restrict this character to get a character
of
    $$
    \{u \in X_*\tq \sum_{\tau\in\op{Gal}(F_1/F)} \tau u = 0\}.
    $$
It can be shown that the character $s$ is trivial on
    $$
    \{u \in X_* \tq \exists \tau\in \op{Gal}(F_1/F)~~ \exists v \in
    X_*.~~ u = \tau v - v\}.
    $$
Thus, by the Tate-Nakayama isomorphism, the character $s$
determines a character $\kappa$ of the cohomology group
    $$
    H^1(\op{Gal}(\bar F/F),T).
    $$

In this way, each cocycle $a_\tau$ gives a complex constant
$\kappa(a_\tau)\in\ring{C}^\times$.

\begin{example}  The element $s\in \ring{C}^\times$ giving the
endoscopic group $H=U_E(1)$ of $SL(2)$ is $s=-1$, which may be
identified with the character $n\mapsto (-1)^n$ of $\ring{Z}$.
This gives the nontrivial character $\kappa$ of
    $$
    H^1(\op{Gal}(\bar F/F),U_E(1)) \cong \ring{Z}/2\ring{Z}.
    $$
\end{example}

\section{Statement of the Fundamental Lemma}

\subsection{Context}\label{sec:context}

Let $G$ be an unramified connected reductive group over $F$.   Let
$H$ be an unramified endoscopic group of $G$. Let $\gamma\in H(F)$
be a strongly regular semisimple element. Let $T_H = C_H(\gamma)$,
and let $T_G$ be a Cartan subgroup of $G$ that is isomorphic to
it.  More details will be given below about how to choose $T_G$.
The choice of $T_G$ matters! Let $\gamma\in T_H(F)$ map to
$\gamma_0\in T_G(F)$ under this isomorphism.

By construction, $\gamma_0$ is semisimple. However, as $G$ may
have more roots than $H$, it is possible for $\gamma_0$ to be
singular, even when $\gamma$ is strongly regular.  If $\gamma\in
H(F)$ is a strongly regular semisimple element with the property
that $\gamma_0$ is also strongly regular, then we will call
$\gamma$ a {\it strongly $G$-regular} element of $H(F)$.

If $\gamma'$ is stably conjugate to $\gamma_0$ with cocycle
$a_\tau$, then $s\in\op{Hom}(X_*,\ring{C}^\times)$ gives
$\kappa(a_\tau)\in \ring{C}^\times$.

Let $K_G$ and $K_H$ be hyperspecial maximal compact subgroups of
$G$ and $H$.  Let $\chi_{G,K}$ and $\chi_{H,K}$ be the
characteristic functions of these hyperspecial subgroups. Set
    \begin{equation}
    \Lambda_{G,H}(\gamma) =
    \left(\prod_{\alpha\in\Phi_G}
    |\alpha(\gamma_0)-1|^{1/2}\right)
    \left[{\frac{\op{vol}(K_T,dt)}{\op{vol}(K,dg)}}\right]
     \sum_{\gamma'\sim\gamma_0}
     \kappa(a_\tau)
    \int_{C_G(\gamma',F)\backslash G(F)}\chi_{G,K} (g^{-1}\gamma' g)
    {\frac{dg}{dt'}}.
    \label{eqn:kappa}
    \end{equation}
The set of roots $\Phi_G$ are taken to be those relative to $T_G$.
The sum runs over all stable conjugates $\gamma'$ of $\gamma_0$,
up to conjugacy. This is a finite sum. The group $K_T$ is defined
to be the maximal compact subgroup of $T_G$. Equation
\ref{eqn:kappa} is a finite linear combination of orbital
integrals (that is, integrals over conjugacy classes in the group
with respect to an invariant measure).   The Haar measures $dt'$
on $C_G(\gamma',F)$ and $dt$ on $T_G(F)$ are chosen so that stable
conjugacy between the two groups is measure preserving.  This
particular linear combination of integrals is called a
$\kappa$-orbital integral because of the term $\kappa(a_\tau)$
that gives the coefficients of the linear combination.  Note that
the integration takes place in the group $G$, and yet the
parameter $\gamma$ is an element of $H(F)$.

The volume terms $\op{vol}(K,dg)$ and $\op{vol}(K_T,dt)$ serve no
purpose other than to make the entire expression independent of
the choice of Haar measures $dg$ and $dt$, which are only defined
up to a scalar multiple.

We can form an analogous linear combination of orbital integrals
on the group $H$.  Set
    \begin{equation}
    \Lambda^{st}_{H}(\gamma) =
    \left(\prod_{\alpha\in\Phi_H} |\alpha(\gamma)-1|^{1/2}\right)
    \left[{\frac{\op{vol}(K_T,dt)}{\op{vol}(K_H,dh)}}\right]
    \sum_{\gamma'\sim\gamma}
    \int_{C_H(\gamma',F)\backslash H(F)}\chi_{H,K} (h^{-1}\gamma' h)
     {\frac{dh}{dt'}}.
    \label{eqn:stable}
    \end{equation}
This linear combination of integrals is like
$\Lambda_{G,H}(\gamma)$, except that $H$ replaces $G$, $K_H$
replaces $K_G$, $\Phi_H$ (taken relative to $T_H$) replaces
$\Phi_G$, and so forth. Also, the factor $\kappa(a_\tau)$ has been
dropped. The linear combination of Equation~\ref{eqn:stable} is
called a stable orbital integral, because it extends over all
stable conjugates of the element $\gamma$ without the factor
$\kappa$. The superscript {\it st\/} in the notation is for
`stable.'

\begin{conjecture} (The fundamental lemma)  For every $\gamma\in
H(F)$ that is strongly $G$-regular semisimple,
    $$
    \Lambda_{G,H}(\gamma) = \Lambda^{st}_H(\gamma).
    $$
\end{conjecture}

\begin{remark} There have been serious efforts over the past twenty
years to prove the fundamental lemma.  These efforts have not yet
led to a proof.  Thus, the fundamental lemma is not a lemma; it is
a conjecture with a misleading name.  Its name leads one to
speculate that the authors of the conjecture may have severely
underestimated the difficulty of the conjecture.
\end{remark}

\begin{remark} Special cases of the fundamental lemma have been
proved.  The case $G=SL(n)$ was proved by Waldspurger \cite{Wal}.
Building on the work of \cite{GKM1}, Laumon has proved that the
fundamental lemma for $G=U(n)$ follows from a purity conjecture
\cite{Lau1}.  The fundamental lemma has not been proved for any
other general families of groups.  The fundamental lemma has been
proved for some groups $G$ of small rank, such as $SU(3)$ and
$Sp(4)$. See \cite{BR}, \cite{Ha1}, \cite{Ha4}.
\end{remark}

\subsection{The significance of the fundamental lemma}

The Langlands program predicts correspondences
$\pi\leftrightarrow\pi'$  between the representation theory of
different reductive groups.  There is a local program for the
representation theory of reductive groups over locally compact
fields, and a global program for {\it automorphic} representations
of reductive groups over the adele rings of global fields.

The Arthur-Selberg trace formula has emerged as a powerful tool in
the Langlands program.  In crude terms, one side of the trace
formula contains terms related to the characters of automorphic
representations.  The other side contains terms such as orbital
integrals.   {\it Thanks to the trace formula, identities between
orbital integrals on different groups imply identities between the
representations of the two groups.}

It is possible to work backwards: from an analysis of the terms in
the trace formula and a precise conjecture in representation
theory, it is possible to make precise conjectures about
identities of orbital integrals. The most basic identity that
appears in this way is the fundamental lemma, articulated above.

The proofs of many major theorems in automorphic representation
theory depend in one way or another on the proof of a fundamental
lemma.  For example, the proof of Fermat's Last Theorem depends on
Base Change for $GL(2)$, which in turn depends on the fundamental
lemma for cyclic base change \cite{Lan1}.  The proof of the local
Langlands conjecture for $GL(n)$ depends on automorphic induction,
which in turn depends on the fundamental lemma for $SL(n)$
\cite{HT}, \cite{HH}, \cite{Wal}.  Properties of the zeta function
of Picard modular varieties depend on the fundamental lemma for
$U(3)$ \cite{R}, \cite{BR}.  Normally, the dependence of a major
theorem on a particular lemma would not be noteworthy. It is only
because the fundamental lemma has not been proved in general, and
because the lack of proof has become a serious impediment to
progress in the field, that the conjecture has become the subject
of increased scrutiny.

\section{Reductions}

To give a trivial example of the fundamental lemma, if $\gamma$
and $\gamma_0$ and their stable conjugates are not in any compact
subgroup, then
    $$\chi_{G,K}(g^{-1}\gamma' g) = 0 \text{ and }
    \chi_{H,K}(h^{-1}\gamma' h) = 0
    $$
so that both $\Lambda_{G,H}(\gamma)$ and
$\Lambda_{H}^{st}(\gamma)$ are zero.  Thus, the fundamental lemma
holds for trivial reasons for such $\gamma$.

\subsection{Topological Jordan decomposition}

A somewhat less trivial reduction of the problem is provided by
the topological Jordan decomposition.  Suppose that $\gamma$ lies
in a compact subgroup.  It can be written uniquely as a product
    $$
    \gamma = \gamma_s \gamma_u = \gamma_u\gamma_s,
    $$
where $\gamma_s$ has finite order, of order prime to the residue
field characteristic $p$, and $\gamma_u$ is topologically
unipotent.  That is,
    $$
    \lim_{n\to\infty} \gamma_u^{p^n} = 1.
    $$
The limit is with respect to the $p$-adic topology.  A special
case of the topological Jordan decomposition $\gamma\in
O_F^\times\subset \ring{G}_m(F)$ is treated in \cite[p20]{Iw}.  In
that case, $\gamma_s$ is defined by the formula
    $$
    \gamma_s = \lim_{n\to\infty} \gamma^{q^n}.
    $$

Let $\gamma$, $\gamma_0$, and $\gamma'$ be chosen as in
Section~\ref{sec:context}.  Each of these elements has a
topological Jordan decomposition.  Let $G_s = C_G(\gamma_{0 s})$
and $H_s = C_H(\gamma_s)$. It turns out that $G_s$ is an
unramified reductive group with unramified endoscopic group $H_s$.
Descent for orbital integrals gives the formulas \cite{LSD}
\cite{Ha2}
    $$
    \begin{array}{lll}
    \Lambda_{G,H}(\gamma) &= \Lambda_{G_s,H_s}(\gamma_u)\\
    &\\
    \Lambda_{H}^{st}(\gamma) &=\Lambda_{H_s}^{st}(\gamma_u).
    \end{array}
    $$
This reduces the fundamental lemma to the case that $\gamma$ is a
topologically unipotent elements.

\subsection{Lie algebras}\label{sec:lie}

It is known (at least when the $p$-adic field $F$ has
characteristic zero), that the fundamental lemma holds for fields
of arbitrary residual characteristic provided that it holds when
the $p$-adic field has sufficiently large residual characteristic
\cite{Ha3}.  Thus, if we are willing to restrict our attention to
fields of characteristic zero, we may assume that the residual
characteristic of $F$ is large.  In fact, in our discussion of a
reduction to Lie algebras in this section, we simply assume that
the characteristic of $F$ is zero.

A second reduction is based on Waldspurger's homogeneity results
for classical groups.  (Homogeneity results have since been
reworked and extended to arbitrary reductive groups by DeBacker,
again assuming mild restrictions on $G$ and $F$.)

When the residual characteristic is sufficiently large, there is
an exponential map from the Lie algebra to the group that has
every topologically unipotent element in its image.   Write
    $$
    \gamma_u = \exp( X),
    $$
for some element $X$ in the Lie algebra.  We may then consider the
behavior of orbital integrals along the curve $\exp(\lambda^2 X)$.
A difficult result of Waldspurger for classical groups states that
if $|\lambda|\le1$, then
    $$
    \begin{array}{lll}
    \Lambda_{G,H}(\exp(\lambda^2 X)) &= \sum a_i |\lambda|^i\\
    \Lambda_{H}^{st}(\exp(\lambda^2 X))&=\sum b_i |\lambda|^i;
    \end{array}
    $$
that is, both sides of the fundamental lemma identity are
polynomials in $|\lambda|$.  If a polynomial identity holds when
$|\lambda|<\epsilon$ for some $\epsilon>0$, then it holds for all
$|\lambda|\le1$.  In particular, it holds at $\gamma_u$ for
$\lambda=1$.  The polynomial growth of orbital integrals makes it
possible to prove the fundamental lemma in a small neighborhood of
the identity element, and then conclude that it holds in general.
In this manner, the fundamental lemma can be reduced to a
conjectural identity in the Lie algebra.

\section{The problem of base points}

The fundamental lemma was formulated above with one omission: we
never made precise how to fix an isomorphism $T_H\leftrightarrow
T_G$ between Cartan subgroups in $H$ and $G$.  Such isomorphisms
exist, because the two Cartan subgroups have the same root data.
But the statement of the fundamental lemma is sensitive to how an
isomorphism is selected between $T_H$ and a Cartan subgroup of
$G$.  If we change the isomorphism, we change the $\kappa$-orbital
integral by a root of unity $\zeta\in\ring{C}^\times$. The
correctly chosen isomorphism will depend on the element $\gamma\in
H(F)$.

The ambiguity of isomorphism was removed by Langlands and Shelstad
in \cite{LS}.  They define a {\it transfer factor\/}
$\Delta(\gamma_H,\gamma_G)$, which is a complex valued function on
$H(F)\times G(F)$. The transfer factor can be defined to have the
property that it is zero unless $\gamma_H\in H(F)$ is strongly
regular semisimple, $\gamma_G\in G(F)$ is strongly regular
semisimple, and there exists an isomorphism (preserving character
groups) from the centralizer of $\gamma_H$ to the centralizer of
$\gamma_G$. There exists $\gamma_0\in G(F)$ such that
    \begin{equation}
    \Delta(\gamma_H,\gamma_0)=1.
    \label{eqn:delta}
    \end{equation}
The correct formulation of the fundamental lemma is to pick the
base point $\gamma_0\in G(F)$ so that Condition~\ref{eqn:delta}
holds.

For classical groups, Waldspurger gives a simplified formula for
the transfer factor $\Delta$ in \cite{WA}.  Furthermore, because
of the reduction of the fundamental lemma to the Lie algebra
(Section~\ref{sec:lie}), the transfer factor may be expressed as a
function on the Lie algebras of $G$ and $H$, rather than as a
function on the group.

\subsection{Base points for unitary groups}\label{sec:bp}

More recently, Laumon (while working on the fundamental lemma for
unitary groups) observed a similarity between Waldspurger's
simplified formula for the transfer factor and the explicit
formula for differents that is found in \cite{SLF}. In this way,
Laumon found a simple description of the matching condition
$\gamma\leftrightarrow\gamma_0$ implicit in the statement of the
fundamental lemma.

\section{Geometric Reformulations of the Fundamental Lemma}
\label{sec:geometry}

From early on, those trying to prove the fundamental lemma have
sought geometric interpretations of the identities of orbital
integrals.  Initially these geometric interpretations were rather
crude.  In the hands of Goresky, Kottwitz, MacPherson, and Laumon
these geometric interpretations have become increasingly
sophisticated. \cite{GKM1}, \cite{GKM2}, \cite{Lau1}, \cite{Lau2}.

This paper is intended to give an introduction to the fundamental
lemma, and the papers giving a geometric interpretation of the
fundamental lemma do not qualify as introductory material.  In
this section, we will be content to describe the geometric
interpretation in broad terms.

\subsection{Old-style geometric interpretations: buildings}

We begin with a geometric interpretation of the fundamental lemma
that was popular in the late seventies and early eighties.  It was
eventually discarded in favor of other approaches when the
combinatorial difficulties became too great.

This approach is to use the geometry of the Bruhat-Tits building
to understand orbital integrals.  We illustrate the approach with
the group $G=SL(2)$.  The term $\chi_{G,K}(g^{-1}\gamma' g)$ that
appears in the fundamental lemma can be manipulated as follows:
    $$
    \begin{array}{lll}
    \chi_{G,K}(g^{-1}\gamma' g)\ne 0 &\Leftrightarrow g^{-1}\gamma' g\in
    K\\
    &\Leftrightarrow \gamma' g\in g K\\
    &\Leftrightarrow \gamma' (g K) = (g K)\\
    &\Leftrightarrow g K \text{ is a fixed point of } \gamma'
    \text{ on } G(F)/K.
    \end{array}
    $$
The set $G(F)/K$ is in bijective correspondence with a set of
vertices in the Bruhat-Tits building of $SL(2)$.  Thus, we may
interpret the orbital integral geometrically as the number of
fixed points of $\gamma'$ in the building that are vertices of a
given type.

Under this interpretation, it is possible to use counting
arguments to obtain explicit formulas for orbital integrals as a
function of $\gamma'$.  In this way, the fundamental lemma was
directly verified for a few groups of small rank such as $SL(2)$
and $U(3)$.

\subsection{Affine grassmannians}

Until the end of Section~\ref{sec:geometry}, let $F = k((t))$, a
field of formal Laurent series.  Except for the discussion of the
results of Kazhdan and Lusztig, the field $k$ will be taken to be
a finite field: $k=\ring{F}_q$.

In 1988, Kazhdan and Lusztig showed that if $F=\ring{C}((t))$,
then $G(F)/K$ can be identified with the points of an ind-scheme
(that is, an inductive limit of schemes) \cite{KL}. This
ind-scheme is called the affine Grassmannian.  The set of fixed
points of an element $\gamma$ can be identified with the set of
points of a scheme over $\ring{C}$, known as the affine Springer
fiber.  The corresponding construction over $\ring{F}_q((t))$ is
mentioned briefly in the final paragraphs of their paper. Rather
than counting fixed points in the building, orbital integral can
be computed by counting the number of points on a scheme over
$\ring{F}_q$.

Based on a description of orbital integrals as the number of
points on schemes over finite fields,  Kottwitz, Goresky, and
MacPherson give a geometrical formulation of the fundamental
lemma.  Furthermore, by making a thorough investigation of the
equivariant cohomology of these schemes, they prove the
geometrical conjecture when $\gamma$ comes from an unramified
Cartan subgroup \cite{GKM1}.

\subsection{Geometric interpretations}

Each of the terms in the fundamental lemma has a nice geometric
interpretation.    Let us give a brief description of the
geometrical counterpart of each term in the fundamental lemma.  We
work with the unitary group, so that we may include various
insights of Laumon.

The geometrical counterpart of cosets $gK$  are self-dual lattices
in a vector space $V$ over $F$.

The counterpart of the support set, $\op{SUP}=\{g\tq g^{-1}\gamma
g \in K\}$, is the affine Springer fiber $X_\gamma$.

The counterpart of the integral of the support set $\op{SUP}$ over
$G$ is counting points on the scheme $X_\gamma$.  The integral
over all of $G$ diverges and the number of fixed points on the
scheme is infinite. For that reason the orbital integral is an
integral over $T\backslash G$, where $T$ is the centralizer of
$\gamma$, rather than over all of $G$.

The counterpart of the integral over $T\backslash G$ is counting
points on a quotient space $Z_\gamma = X_\gamma/\ring{Z}^\ell$.
(There is a free action of a group $\ring{Z}^\ell$ on $X_\gamma$,
and $Z_\gamma$ is the quotient.)

The geometric counterpart of $\kappa(a_\tau)$ is somewhat more
involved. For elliptic endoscopic groups of unitary groups
$\kappa$ has order $2$.  The character $\kappa$ has the form.
    $$
    \kappa:H^1(\op{Gal}(\bar F/F),T)\cong
    (\ring{Z}/2\ring{Z})^\ell \to \{\pm1\}.
    $$
The character $\kappa$ pulls back to a character of
$\ring{Z}^\ell$. The rational points of $X_\gamma$ are identified
with self-dual lattices: $A^{\perp} = A$. The points of the
quotient space $Z_\gamma$ are lattices that are self-dual modulo
the group action: $A^{\perp} =\lambda\cdot A$, for some
$\lambda\in \ring{Z}^\ell$. The character $\kappa$ then partitions
the points of $Z_\gamma$ into two sets, depending on the sign of
$\kappa(\lambda)$:
    $$
    Z_\gamma^{\pm} = \{A ~|~ A^{\perp} = \lambda A;\quad
    \kappa(\lambda) = \pm 1\}.
    $$
(In a more sophisticated treatment of $\kappa(a_\tau)$, it gives
rise to a local system on $Z_\gamma$; and counting points on
varieties gives way to Grothendieck's trace formula.)

The counterpart of the $\kappa$-orbital integral
$\Lambda_{G,H}(\gamma)$ is the number
    $$
    \#Z_\gamma^+ - \#Z_\gamma^-.
    $$

The counterpart of the stable-orbital integral
$\Lambda_H^{st}(\gamma)$ is the number
    $$
    \#Z_\gamma^{H,st}
    $$
for a corresponding variety constructed from the endoscopic group.

The factors $\prod_{\Phi} |\alpha(\gamma)-1|^{1/2}$ that appear on
the two sides of the fundamental lemma can be combined into a
single term
    $$
    \prod_{\alpha\in\Phi_G\setminus
    \Phi_H}|\alpha(\gamma)-1|^{1/2}.
    $$
This has the form $q^{-d}$ for some value $d=d(\gamma)$.  The
factor $q^{-d}$ has been interpreted in various ways.  We mention
that \cite{LR} interprets $q^d$ as the points on an affine space
of dimension $d$.  That paper expresses the hope that it might be
possible to find an embedding $Z_\gamma^-\to Z_\gamma^+$ such that
the complement of the embedded $Z_\gamma^-$ in $Z_\gamma^+$ is a
rank $d$ fiber bundle over $Z_\gamma^{H,st}$. The realization of
this hope would give an entirely geometric interpretation of the
fundamental lemma.  Laumon and Rapoport found that this
construction works over $\ring{F}_{q^2}((t))$, but not over
$\ring{F}_q((t))$. In more recent work of Laumon, the constant $d$
is interpreted geometrically as the intersection multiplicity of
two singular curves.

\subsection{Compactified Jacobians}

Laumon, in the case of unitary groups, has made the splendid
discovery that the orbital integrals -- as they appear in the
fundamental lemma -- count points on the compactification of the
Jacobians of a singular curve associated with the semisimple
element $\gamma$. (In fact, $Z_\gamma$ is homeomorphic to and can
be replaced with the compactification of a Jacobian.)  Thus, the
fundamental lemma may be reformulated as a relation between the
compactified Jacobians of these curves.  By showing that the
singular curve for the endoscopic group $H$ is a perturbation of
the singular curve for the group $G$, he is able relate the
compactified Jacobians of the two curves, and prove the
fundamental lemma for unitary groups (assuming a purity hypothesis
related to the cohomology of the schemes).

\begin{figure}[htb]
  \centering
  \includegraphics{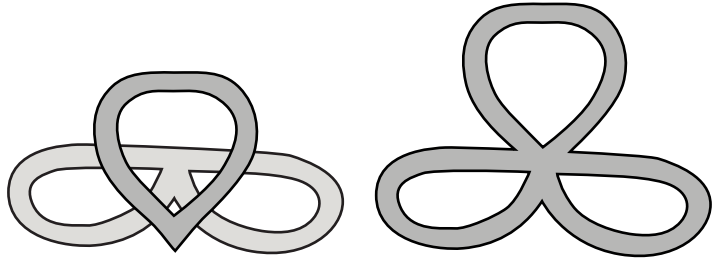}
  \caption{The singular curve on the left can be deformed into
  the singular curve on the right by pulling up on the center
  ring. The curve on the left controls $\Lambda^{st}_H(\gamma)$,
  and the curve on the right controls $\Lambda_{G,H}(\gamma)$.
  This deformation relating the two curves
  is a key part of Laumon's work on the fundamental
  lemma for unitary groups.}
  \label{fig:gfcchcp}
\end{figure}

The origin of the curve $C$ is the following.  The ring
$O_F[\gamma]$ is the completion at a point of the local ring of a
curve $C$.  In the interpretation in terms of Jacobians, the
self-dual lattices $A^{\perp} = A$ that appear in the geometric
interpretation above are replaced with $O_C$-modules, where $O_C$
is the structure sheaf of $C$.

The audio recording of Laumon's lecture at the Fields Institute on
this research is highly recommended \cite{Lau3}.

\subsection{Final remarks}

\begin{remark}
The fundamental lemma is an open ended problem, in the sense that
as researchers develop new trace formulas (the symmetric space
trace formula \cite{Jac}, the twisted trace formula \cite{KS}, and
so forth) and as they compare trace formulas for different groups,
it will be necessary to formulate and prove generalized versions
of the fundamental lemma.  The version of the fundamental lemma
stated in this paper should be viewed as a template that should be
adapted according to an evolving context.
\end{remark}

\begin{remark}
The methods of Goresky, Kottwitz, MacPherson, and Laumon are
limited to fields of positive characteristic.  This may at first
seem to be a limitation of their method.  However, there are ideas
about how to use motivic integration to lift their results from
positive characteristic to characteristic zero (see \cite{CH}).
Waldspurger also has results about lifting to characteristic zero
that were presented at the Labesse conference, but I have not seen
a preprint \cite{CL}.
\end{remark}

\begin{remark}
In some cases, it is now known how to deduce stronger forms of the
fundamental lemma from weaker versions.  For example, it is known
how to go from the characteristic function of the hyperspecial
maximal compact groups to the full Hecke algebra \cite{Ha3}.   A
descent argument replaces twisted orbital integrals by ordinary
orbital integrals.   However, relations between weighted orbital
integrals remain a serious challenge.
\end{remark}

\begin{remark}
There has been much research on the fundamental lemma that has not
been discussed in detail in this paper, including other forms of
the fundamental lemma. For just one example, see \cite{Ngo} for
the fundamental lemma of Jacquet and Ye.  Other helpful references
include \cite{LD} and \cite{WICM}.
\end{remark}

\bibliographystyle{amsplain}

    \renewcommand{\thefootnote}{}
    \footnote{This work was supported in part by the NSF.}
    \footnote{This paper is based on lectures at the
    Fields Institute on June 25--27, 2003.  Slides and audio
    are available at\hfil\break\qquad
    {\tt http://www.fields.utoronto.ca/audio/02-03/\#CMI\_summer\_school}}.

    \footnote{Copyright 2003, Thomas C. Hales.}
    \footnote{
This work is licensed under the Creative Commons Attribution
License. To view a copy of this license, visit
http://creativecommons.org/licenses/by/1.0/ or send a letter to
Creative Commons, 559 Nathan Abbott Way, Stanford, California
94305, USA.}

\end{document}